\def\DateTime{11/December/2006, 15:00(JP)}
\def\Version{Version $1.0$}

\def\no{\if01}
\def\iftwelvept{\no}

\def\ifusepdf{\no}
\def\ifpsfont{\no}

\iftwelvept
\documentclass[leqno,12pt]{amsart}
\else
\documentclass[leqno]{amsart}
\fi
\usepackage{amssymb}
\usepackage{amscd}
\usepackage{latexsym}
\usepackage{verbatim}
\usepackage[all]{xy}

\ifusepdf
\usepackage{hyperref}
\else\fi
\ifpsfont
\usepackage[T1]{fontenc}
\usepackage{times}
\else\fi

\iftwelvept
\else\fi

\theoremstyle{plain}
\newtheorem{Theorem}{Theorem}[section]
\newtheorem{Proposition}[Theorem]{Proposition}
\newtheorem{Lemma}[Theorem]{Lemma}
\newtheorem{Corollary}[Theorem]{Corollary}

\newtheorem{Claim}{Claim}[Theorem]

\theoremstyle{definition}

\renewcommand{\theTheorem}{\arabic{section}.\arabic{Theorem}}
\renewcommand{\theClaim}{\arabic{section}.\arabic{Theorem}.\arabic{Claim}}
\renewcommand{\theequation}{\arabic{section}.\arabic{Theorem}.\arabic{Claim}}

\def\rom{\textup}
\newcommand{\ZZ}{{\mathbb{Z}}}
\newcommand{\QQ}{{\mathbb{Q}}}

\newcommand{\CC}{{\mathbb{C}}}

\newcommand{\OO}{{\mathcal{O}}}
\newcommand{\Proj}{\operatorname{Proj}}

\newcommand{\Ker}{\operatorname{Ker}}

\newcommand{\codim}{\operatorname{codim}}

\newcommand{\rank}{\operatorname{rk}}
\newcommand{\length}{\operatorname{length}}
\newcommand{\ord}{\operatorname{ord}}
\newcommand{\acherncl}{\widehat{{c}}}
\newcommand{\aChow}{\widehat{\operatorname{CH}}}
\newcommand{\adeg}{\widehat{\operatorname{deg}}}

\newcommand{\sgn}{\operatorname{sgn}}

\newcommand{\Sym}{\operatorname{Sym}}

\newcommand{\Proof}{{\sl Proof.}\quad}
\newcommand{\QED}{{\unskip\nobreak\hfil\penalty50\quad\null\nobreak\hfil
{$\Box$}\parfillskip0pt\finalhyphendemerits0\par\medskip}}
\newcommand{\rest}[2]{\left.{#1}\right\vert_{{#2}}}

\begin{document}

\title[Subsheaves of a hermitian torsion free coherent sheaf]%
{Subsheaves of a hermitian torsion free
coherent sheaf on an arithmetic variety}
\author{Atsushi Moriwaki}
\address{Department of Mathematics, Faculty of Science,
Kyoto University, Kyoto, 606-8502, Japan}
\email{moriwaki@math.kyoto-u.ac.jp}
\date{\DateTime, (\Version)}


\maketitle


\section*{Introduction}
\renewcommand{\theTheorem}{\Alph{Theorem}}

Let $K$ be a number field and $O_K$ the ring of integers of $K$.
Let $(E, h)$ be a hermitian finitely generated flat $O_K$-module.
For an $O_K$-submodule $F$ of $E$,
let us denote by $h_{F\hookrightarrow E}$ the submetric of $F$ induced by $h$.
It is well known that the set of all saturated $O_K$-submodules
$F$ with $\adeg(F, h_{F \hookrightarrow E}) \geq c$ is finite for any real numbers $c$
(for details, see \cite[the proof of Proposition~3.5]{MoBo}).

In this note, we would like to give its generalization on a projective arithmetic variety.
Let $X$ be a normal and projective arithmetic variety.
Here we assume that $X$ is an arithmetic surface
to avoid several complicated technical definitions on a higher dimensional
arithmetic variety.
Let us fix  a nef and 
big $C^{\infty}$-hermitian
invertible sheaf $\overline{H}$ on $X$ as a polarization of $X$.
Then we have the following finiteness of saturated subsheaves with bounded arithmetic degree,
which is also a generalization of a partial result 
\cite[Corollary~2.2]{MoUBo}.

\begin{Theorem}[cf. Theorem~\ref{thm:finite:subsheaf}]
\label{thm:finite:subsheaf:intro}
Let $E$ be a torsion free coherent sheaf on $X$ and 
$h$ a $C^{\infty}$-hermitian metric of $E$ on $X(\CC)$.
For any real number $c$, the set of all saturated $\OO_{X}$-subsheaves $F$ of $E$
with $\adeg( \acherncl_1(\overline{H}) \cdot \acherncl_1(F, h_{F\hookrightarrow E})) \geq c$ 
is finite.
\end{Theorem}

For a non-zero $C^{\infty}$-hermitian torsion free coherent sheaf $\overline{G}$ on $X$,
the {\em arithmetic slope} $\hat{\mu}_{\overline{H}}(\overline{G})$ of $\overline{G}$
with respect to $\overline{H}$ is defined by
\[
\hat{\mu}_{\overline{H}}(\overline{G}) =
\frac{\adeg (\acherncl_1(\overline{H}) \cdot \acherncl_1(\overline{G}) )}{\rank G}.
\]
As defined in the paper \cite{MoUBo},
$(E, h)$ is said to be {\em arithmetically $\mu$-semistable} with respect to $\overline{H}$ if,
for any non-zero saturated $\OO_{X}$-subsheaf $F$ of $E$,
\[
\hat{\mu}_{\overline{H}}( F, h_{F\hookrightarrow E})
\leq 
\hat{\mu}_{\overline{H}}(E, h).
\]
The above semistability yields
an arithmetic analogue of the Harder-Narasimham filtration of
a torsion free sheaf on an algebraic variety as follows:
A filtration
\[
0 = E_0 \subsetneq E_1 \subsetneq \cdots \subsetneq E_l = E
\]
of $E$ is called an {\em arithmetic Harder-Narasimham filtration of $(E, h)$
with respect to $\overline{H}$} if
\begin{enumerate}
\renewcommand{\labelenumi}{(\arabic{enumi})}
\item
$E_i/E_{i-1}$ is torsion free for every $1 \leq i \leq l$.

\item
Let $h_{E_i/E_{i-1}}$ be a $C^{\infty}$-hermitian metric of $E_i/E_{i-1}$
induced by $h$, that is,
\[
h_{E_i/E_{i-1}} = (h_{E_i \hookrightarrow E})_{E_i \twoheadrightarrow E_i/E_{i-1}} = 
(h_{E \twoheadrightarrow E/E_{i-1}})_{E_i/E_{i-1}  \hookrightarrow E/E_{i-1}}
\]
(for details, see Proposition~\ref{prop:sub:quot:metric}).
Then $(E_i/E_{i-1}, h_{E_i/E_{i-1}})$ is arithmetically $\mu$-semistable with respect to $\overline{H}$.

\item
$\hat{\mu}_{\overline{H}}(E_1/E_0, h_{E_1/E_0}) >
\hat{\mu}_{\overline{H}}(E_2/E_1, h_{E_2/E_1}) > \cdots >
\hat{\mu}_{\overline{H}}(E_l/E_{l-1}, h_{E_l/E_{l-1}})$.
\end{enumerate}
As a consequence of the above theorem, we can show
the unique existence of an arithmetic Harder-Narasimham filtration:

\begin{Theorem}[cf. Theorem~\ref{thm:existence:HN}]
\label{thm:HN:filt:intro}
There is a unique arithmetic Harder-Narasimham filtration of $(E, h)$.
\end{Theorem}

\renewcommand{\theTheorem}{\arabic{section}.\arabic{subsection}.\arabic{Theorem}}
\renewcommand{\theClaim}{\arabic{section}.\arabic{subsection}.\arabic{Theorem}.\arabic{Claim}}
\renewcommand{\theequation}{\arabic{section}.\arabic{subsection}.\arabic{Theorem}.\arabic{Claim}}

\section{Preliminaries}
\subsection{Hermitian vector space}
\setcounter{Theorem}{0}
In this subsection, let us recall several basic facts of hermitian complex vector spaces.

Let $(V, h)$ be a finite dimensional hermitian complex vector space, i.e.,
$V$ is a finite dimensional vector space over $\CC$ and $h$ is
a hermitian metric of $V$.
Let $\phi : V' \to V$ be an injective homomorphism of complex vector spaces.
If we set $h'(x, y) = h(\phi(x), \phi(y))$, then
$h'$ is a hermitian metric of $V'$.
This metric $h'$ is called the {\em submetric of $V'$ induced by $h$ and $V' \to V$},
and it is denoted by $h_{V' \hookrightarrow V}$.

Let $\psi : V \to V''$ be a surjective homomorphism of complex vector spaces.
Let $W$ be the orthogonal complement of $\Ker(\psi)$ with respect to
$h$. Let $h_{W \hookrightarrow V}$ be the submetric of $W$ induced by
$h$ and $W \to V$.
Then there is a unique hermitian metric $h''$ of $V''$ such that
the isomorphism $\rest{\psi}{W} : W \to V''$ gives rise to
an isometry $(W, h_{W \hookrightarrow V}) \overset{\sim}{\longrightarrow} (V'', h'')$.
The metric $h''$ is called  the {\em quotient metric of $V''$ induced by $h$ and $V \to V''$},
and it is denoted by $h_{V \twoheadrightarrow V''}$.

For simplicity, the submetric $h_{V' \hookrightarrow V}$ and
the quotient metric $h_{V \twoheadrightarrow V''}$ are often denoted by
$h_{V'}$ and $h_{V''}$ respectively.
It is easy to see the following proposition:

\begin{Proposition}
\label{prop:sub:quot:metric}
Let $V, V', V''$ be finite dimensional complex vector spaces with
$V'' \subseteq V' \subseteq V$. Let $h$ be a hermitian metric of $V$.
Then
\[
(h_{V' \hookrightarrow V})_{V' \twoheadrightarrow V'/V''} =
(h_{V \twoheadrightarrow V/V''})_{V'/V'' \hookrightarrow V/V''}
\]
as hermitian metrics of $V'/V''$.
\end{Proposition}

More generally, we have the following lemma:

\begin{Lemma}
\label{lem:comp:inq:two:metrics}
Let $(V, h)$ be a finite dimensional hermitian complex vector space.
Let $W$ and $U$ be subspaces of $V$.
Let us consider a natural homomorphism
\[
\phi : W \hookrightarrow V \to V/U
\]
of complex vector spaces. Let $Q$ be the image of $\phi$.
Let us consider two natural hermitian metrics $h_1$ and $h_2$ of $Q$ given by
\[
h_1 = (h_{W \hookrightarrow V})_{W \twoheadrightarrow Q}
\quad\text{and}\quad
h_2 = (h_{V \twoheadrightarrow V/U})_{Q \hookrightarrow V/U}.
\]
Then $h_1(x, x) \geq h_2(x, x)$ for all $x \in Q$.
In particular, if $\{ x_1, \ldots, x_s \}$ is a basis of $Q$, then
$\det(h_1(x_i, x_j)) \geq \det(h_2(x_i, x_j))$.
\end{Lemma}

\Proof
Let $T$ be the orthogonal complement of $\Ker(\phi : W \to Q)$
with respect to $h_{W \hookrightarrow V}$.
Then $h(v, v) = h_1(\phi(v), \phi(v))$ for all $v \in T$.
Let $U^{\perp}$ be the orthogonal complement of $U$ with respect to
$h$. Then, for $v \in T$, we can set $v = u + u'$ with
$u \in U$ and $u' \in U^{\perp}$.
Then $h_2(\phi(v), \phi(v)) = h(u', u')$. Thus
\[
h_2(\phi(v), \phi(v)) = h(u', u') \leq h(v, v) = h_1(\phi(v), \phi(v)).
\]
For the last assertion, see \cite[Lemma~3.4]{MoBo}.
\QED

\bigskip
Let $e_1, \ldots, e_n$ be an orthonormal basis of $V$ with respect to $h$.
Let $V^{\vee}$ be the dual space of $V$ and
$e^{\vee}_1, \ldots, e^{\vee}_n$ the dual basis of $e_1, \ldots, e_n$.
For $\phi, \psi \in V^{\vee}$, we set
\[
h^{\vee}(\phi, \psi) = \sum_{i=1}^n a_i \bar{b}_i, 
\]
where $\phi = a_1 e^{\vee}_1 + \cdots + a_n e^{\vee}_n$ and
$\psi = b_1 e^{\vee}_1 + \cdots + b_n e^{\vee}_n$.
It is easy to see that $h^{\vee}$ does not depend on the choice of
the orthonormal basis of $V$, so that
the hermitian metric $h^{\vee}$ of $V^{\vee}$ is called the {\em dual hermitian metric of $h$}.
Moreover we can easily check the following facts:

\begin{Proposition}
\label{prop:dual:hermitian:vector}
\begin{enumerate}
\renewcommand{\labelenumi}{(\arabic{enumi})}
\item
${\displaystyle h^{\vee}(\phi, \phi) = \sup_{ x \in V \setminus \{ 0 \} } \frac{\vert \phi(x) \vert^2}{h(x, x)}}$.

\item
Let $x_1, \ldots, x_n$ be a basis of $V$ and
$x_1^{\vee}, \ldots, x_n^{\vee}$ be the dual basis of $V^{\vee}$.
If we set $H = (h(x_i, x_j))$ and $H^{\vee} = (h^{\vee}(x^{\vee}_i, x^{\vee}_j))$, then
$H^{\vee} = \overline{H}^{-1}$.

\item
Let $0 \to V_1 \to V_2 \to V_3 \to 0$ be an exact sequence of finite dimensional complex vector
spaces and $h_1, h_2, h_3$ hermitian metrics of $V_1, V_2, V_3$ respectively.
We assume that $h_1 = (h_2)_{V_1 \hookrightarrow V_2}$
and $h_3 = (h_2)_{V_2 \twoheadrightarrow V_3}$.
Let us consider the dual exact sequence  $0 \to V^{\vee}_3 \to V^{\vee}_2 \to V^{\vee}_1 \to 0$ of
$0 \to V_1 \to V_2 \to V_3 \to 0$ and the dual hermitian metrics $h^{\vee}_1, h^{\vee}_2, h^{\vee}_3$ of
$h_1, h_2, h_3$ respectively.
Then $h^{\vee}_3 = (h_2^{\vee})_{V^{\vee}_3 \hookrightarrow V^{\vee}_2}$
and $h_1^{\vee} = (h_2^{\vee})_{V^{\vee}_2 \twoheadrightarrow V^{\vee}_1}$.
\end{enumerate}
\end{Proposition}

\bigskip
Let $(U, h_U)$ and $(W, h_W)$ be
finite dimensional hermitian vector spaces over $\CC$.
Then $U \otimes_{\CC} W$ has the standard hermitian metric
$h_{U} \otimes h_{W}$ defined by
\[
(h_{U} \otimes h_{W})(u \otimes w, u' \otimes w') = h_U(u, u') h_W(w, w').
\]
Thus the standard hermitian metric of $\bigotimes^r V$ is given by 
\[
(\bigotimes^r h)(v_1 \otimes \cdots v_r, v'_1 \otimes \cdots \otimes v'_r) =
h(v_1, v'_1) \cdots h(v_r, v'_r).
\]
Let $\pi : \bigotimes^r V \to \bigwedge^r V$ 
be the natural surjective homomorphism and
$\bigwedge^r h$ a hermitian metric of $\bigwedge^r V$ given by
\[
\bigwedge^r h = r! (\bigotimes^r h)_{\bigotimes^r V \twoheadrightarrow \bigwedge^r V}.
\]
Then we have the following:

\begin{Proposition}
\label{prop:hermitian:metric:wedge}
$(\bigwedge^r h)(x_1 \wedge \cdots \wedge x_r, x_1 \wedge \cdots \wedge x_r) =
\det (h(x_i, x_j))$.
\end{Proposition}

\Proof
For $a_1, \ldots, a_r \in V$, we set
\[
\phi(a_1, \ldots, a_r) = \frac{1}{r!} \sum_{\sigma \in S_r} \sgn(\sigma)
a_{\sigma(1)} \otimes \cdots \otimes a_{\sigma(r)}.
\]
Then, by an easy calculation, 
for $\sigma \in S_r$ and $a_1, \ldots, a_r, b_1, \ldots, b_r \in V$,
we can see
\addtocounter{Claim}{1}
\begin{multline}
\label{prop:hermitian:metric:wedge:eqn:1}
(\bigotimes^r h)(a_{\sigma(1)} \otimes \cdots \otimes a_{\sigma(r)},
\phi(b_1, \ldots, b_r)) = \\
\sgn(\sigma) 
(\bigotimes^r h)(a_{1} \otimes \cdots \otimes a_{r},
\phi(b_1, \ldots, b_r)) 
\end{multline}
Note that $\Ker(\pi)$ is generated by elements of type
\[
a_1 \otimes \cdots \otimes a_r,
\]
where $a_i = a_j$ for some $i \not= j$.
Therefore, by \eqref{prop:hermitian:metric:wedge:eqn:1},
$\phi(x_1, \ldots, x_r) \in \Ker(\pi)^{\perp}$ for all $x_1, \ldots, x_r \in V$.
Thus, since 
\[
\pi(\phi(x_1, \ldots, x_r)) = x_1 \wedge \cdots \wedge x_r,
\]
we have
\[
(\bigotimes^r h)_{\bigotimes^r V \twoheadrightarrow \bigwedge^r V}(x_1 \wedge \cdots \wedge x_r,
x_1 \wedge \cdots \wedge x_r) = 
(\bigotimes^r h)(\phi(x_1, \ldots, x_r), \phi(x_1, \ldots, x_r)).
\]
On the other hand, by using \eqref{prop:hermitian:metric:wedge:eqn:1} again,
we can check
\[
(\bigotimes^r h)(\phi(x_1, \ldots, x_r), \phi(x_1, \ldots, x_r)) = \frac{1}{r!} \det(h(x_i, x_j)).
\]
Therefore we get our assertion.
\QED

\subsection{Finitely generated modules over a $1$-dimensional noetherian integral domain}
\setcounter{Theorem}{0}
Let $R$ be a noetherian integral domain with $\dim R = 1$, and
$K$ the quotient field of $R$.
For $a\in R \setminus \{ 0 \}$,
we set $\ord_R(a) = \length_R(R/aR)$, which
yields a homomorphism $\ord_R : R \setminus \{ 0\} \to \ZZ$, that is,
$\ord_R(ab) = \ord_R(a) +\ord_R(b)$
for $a, b \in R \setminus \{ 0\}$.
Thus it extends to a homomorphism on $K^{\times}$ given by
$\ord_R(a/b)= \ord_R(a) - \ord_R(b)$.

\begin{Proposition}
\label{prop:base:change:rule}
Let $E$ be a finitely generated $R$-module.
Let $s_1, \ldots, s_r$ and $s'_1, \ldots, s'_r$ be sequences of elements
of $E$ such that $s_1, \ldots, s_r$ and $s'_1, \ldots, s'_r$ form
bases of $E \otimes_R K$ respectively.
Let $A = (a_{ij})$ be an $r \times r$-matrix such that
$a_{ij} \in K$ for all $i, j$ and
$s'_i = \sum_{j=1}^r a_{ij}s_j$ in $E \otimes_R K$ for all $i$.
Then
\[
\length_R(E/R s'_1 + \cdots + R s'_r) = 
\length_R(E/Rs_1 + \cdots + Rs_r)
+ \ord_R(\det(A)).
\]
\end{Proposition}

\Proof
We set $M = Rs_1 + \cdots + Rs_r$ and $M' = Rs'_1 + \cdots + Rs'_r$.
First we assume that $M' \subseteq M$.
Then $a_{ij} \in R$.
An exact sequence
\[
0 \to M/M'
\to E/M'\to E/M \to 0.
\]
yields
\[
\length_R(E/M') =
\length_R(E/M) + \length_R(M/M').
\]
Note that $M$ is a free $R$-module.
Let $\phi : M \to M$ be an endomorphism given by $\phi(s_i) = s'_i$.
Then, by [EGA~IV, Lemme~21.10.17.3],
$\length_R(M/\phi(M)) = \length_R(R/\det(\phi)R)$.
Thus we get
\[
\length_R(E/M') =
\length_R(E/M) + \length_R(R/\det(A)R).
\]

Next we consider a general case.
Since $E/M$ is a torsion module, there is $b \in R \setminus \{ 0\}$
with $bM' \subseteq M$.
Thus, by the previous observation,
\[
\length_R(E/bM') =
\length_R(E/M) + \length_R(R/\det(bA)R)
\]
because $bs_i = \sum_{j=1}^r ba_{ij}s_j$ in $E \otimes_R K$ for all $i$.
Moreover
\[
\length_R(E/bM') =
\length_R(E/M') + \length_R(R/b^rR).
\]
Hence the proposition follows.
\QED

\begin{Corollary}
\label{cor:base:change:rule}
\begin{enumerate}
\renewcommand{\labelenumi}{(\arabic{enumi})}
\item
Let $\{ x_1, \ldots, x_r\}$ be a basis of $E \otimes_R K$.
Let $s_1, \ldots, s_r \in E$ and $a \in R \setminus \{ 0 \}$ such that
$a x_i = s_i$ in $E \otimes_R K$ for all $i$.
Then the number
\[
\length_R(E/Rs_1 + \cdots + Rs_r) - r \ord_R(a)
\]
does not depend on the choice of $s_1, \ldots, s_r$ and $a$, so that it is denoted
by $\ell_R(E; x_1, \ldots, x_r)$.

\item
Let $\{ x_1, \ldots, x_r \}$ and $\{ x'_1, \ldots, x'_r \}$ be bases of
$E \otimes_R K$. Let $B = (b_{ij})$ be an $r \times r$ matrix such that
$x'_i = \sum_{j=1}^r b_{ij} x_j$ for all $i$.
Then
\[
\ell_R(E; x'_1, \ldots, x'_r) = \ell_R(E; x_1, \ldots, x_r) + \ord_R(\det(B)).
\]
\end{enumerate}
\end{Corollary}

\Proof
(1) Let $s'_1, \ldots, s'_r \in E$ and $a' \in R \setminus \{ 0 \}$ be
another choice with
$a' x_i = s'_i$ in $E \otimes_R K$ for all $i$.
Then $s'_i = (a'/a) s_i$ in $E \otimes_R K$.
Thus, by the previous proposition,
\[
\length_R(E/R s'_1 + \cdots + R s'_r) = 
\length_R(E/Rs_1 + \cdots + Rs_r)
+ \ord_R((a'/a)^r),
\]
which yields the assertion.

\medskip
(2) Let us choose $a, b \in R \setminus \{ 0 \}$ and $s_1, \ldots, s_r \in E$ such that
$ax_i = s_i$ in $E \otimes_R K$ for all $i$ and $b b_{ij} \in R$ for all $i, j$.
If we set $s'_i = \sum_j (bb_{ij})s_i$, then $ab x'_i = s'_i$ in $E \otimes_R K$ for all $i$.
Thus
\begin{align*}
\ell_R(E; x_1, \ldots, x_r) & =  \length_R(E/R s_1 + \cdots + R s_r) - r \ord_R(a) \\
\ell_R(E; x'_1, \ldots, x'_r) & =  \length_R(E/R s'_1 + \cdots + R s'_r) - r \ord_R(ab).
\end{align*}
On the other hand, by the previous proposition,
\[
\length_R(E/R s'_1 + \cdots + R s'_r) = \length_R(E/R s_1 + \cdots + R s_r)
+ \ord_R(\det(bB)).
\]
Hence we obtain (2).
\QED

\subsection{Subsheaves of a torsion free coherent sheaf}
\setcounter{Theorem}{0}
In this subsection, we consider how we can get a saturated subsheaf.

\begin{Proposition}
\label{prop:subsheaves:vector:subspaces}
Let $X$ be an irreducible noetherian integral scheme,
$\eta$ the generic point of $X$, and $K = \OO_{X,\eta}$
the function field of $X$.
Let $E$ be a torsion free coherent sheaf on $X$.
Let $\Sigma(X, E)$ be the set of all saturated $\OO_X$-subsheaves
of $E$ and $\Sigma(K, E_{\eta})$ the set of all vector subspaces of $E_{\eta}$ over $K$.
Then the map $\gamma : \Sigma(X, E) \to \Sigma(K, E_{\eta})$
given by $\gamma(F) = F_{\eta}$ is bijective.
For a vector subspace $W$ of $E_{\eta}$ over $K$,
the subsheaf given by $\gamma^{-1}(W)$ is called
the {\em saturated $\OO_X$-subsheaf of $E$ induced by $W$}
and is denoted by $\OO_X(W; E)$.
\end{Proposition}

\Proof
Let us begin with the following lemma:

\begin{Lemma}
\label{unique:saturate:subsheaf}
Let $F, G$ be $\OO_X$-subsheaves of $E$ such that
$F$ is saturated in $E$ and $F_{\eta} = G_{\eta}$.
Then $F \supseteq G$.
\end{Lemma}

\Proof
Let us consider a homomorphism
$\phi : G \to E \to E/F$. Then $\phi_{\eta} = 0$.
Since $E/F$ is torsion free, we have $\phi = 0$,
which means that $G \subseteq F$.
\QED

The injectivity of $\gamma$ is a consequence of the above lemma.
Let $W$ be a vector subspace of $E_{\eta}$ over $K$.
We set $F(U) = W \cap E(U)$ for any Zariski open set $U$ of $X$.
Then $F_{\eta} = W$. We need to see that $F$ is saturated in $E$.
Since $F$ is the kernel of the natural homomorphism
$E \to E_{\eta} \to E_{\eta}/W$, we have an injection
$E/F \hookrightarrow E_{\eta}/W$, so that $E/F$ is torsion free.
\QED

\begin{Proposition}
\label{prop:invertible:section:1to1}
Let $X$ be a noetherian scheme and
$E$ a locally free coherent sheaf on $X$.
Let $\pi : P = \Proj(\bigoplus_{d \geq 0} \Sym^d(E^{\vee})) \to X$
be the projective bundle and
$\OO_P(1)$ the tautological line bundle of $P \to X$.
Let $\Gamma(X, P)$ be the set of all sections of $\pi : P \to X$.
Moreover let $\Sigma'_1(X, E)$ be the set of all
$\OO_X$-subsheaves $L$ 
such that
$L$ is invertible and $E/L$ is locally free.
For $s \in \Gamma(X, P)$,
let 
\[
\phi_s : s^*(\OO_P(-1)) \to s^*\pi^*(E) = E
\]
be a homomorphism
obtained from the dual homomorphism
$\OO_P(-1) \to \pi^*(E)$ of
the natural homomorphism
$\pi^*(E^{\vee}) \to \OO_P(1)$
by applying $s^*$.
We denote the image of
$\phi_s : s^*(\OO_P(-1)) \to E$ by $L(s)$.
Then $L(s) \in \Sigma'_1(X, E)$ for all $s \in \Gamma(X, P)$ and
a map
\[
  \Gamma(X, P) \to \Sigma'_1(X, E)
\]
given by $s \mapsto L(s)$ is bijective.
\end{Proposition}

\Proof
See \cite[Theorem~7.1 and Proposition~7.12]{Harts}.
\QED

\subsection{Hermitian locally free coherent sheaf on a smooth variety}
\setcounter{Theorem}{0}
Let $X$ be a smooth variety over $\CC$,
$\eta$ be the generic point of $X$, and $K =\OO_{X,\eta}$ the function field of $X$.

\begin{Proposition}
\label{prop:coincide:hermitian:locally:free:sheaf}
Let  $(E, h)$ and $(E', h')$ be $C^{\infty}$-hermitian locally free coherent
sheaves on $X$. If there is a dense Zariski open set $U$ of $X$ such that
$\rest{(E, h)}{U}$ is isometric to $\rest{(E', h')}{U}$, then
this isometry extends to an isometry over $X$.
\end{Proposition}

\Proof
Since $V = E_{\eta}$ is isomorphic to $E'_{\eta}$,
we may assume that $E'$ is a subsheaf of $V$.
Then $\rest{(E, h)}{U}$ coincides with $\rest{(E', h')}{U}$.

First let us see that $E = E'$. For this purpose, it is sufficient to see that
$E_{\gamma} = E'_{\gamma}$ for all codimension one points $\gamma$.
Let $\{ \omega_1, \ldots, \omega_r \}$ and
$\{ \omega'_1, \ldots, \omega'_r \}$ be local bases of $E_{\gamma}$ and
$E'_{\gamma}$ respectively.
Then there are $r \times r$-matrices $(a_{ij})$ and $(b_{ij})$ such that
$a_{ij}, b_{ij} \in K$ for all $i, j$ and
\[
\omega'_i = \sum_{j=1}^r a_{ij}\omega_j,\quad
\omega_i = \sum_{j=1}^r b_{ij}\omega'_j
\]
for all $i$.  Clearly $(a_{ij})(b_{ij}) = (b_{ij})(a_{ij}) = (\delta_{ij})$.

\begin{Claim}
$a_{ij}, b_{ij} \in \OO_{X,\gamma}$ for all $i,j$.
\end{Claim}

For each $i$, we set $e_i = \min_{1\leq j \leq r} \{ \ord_{\gamma}(a_{ij}) \}$.
We assume that $e_i < 0$.
Let $t$ be a local parameter of $\OO_{X, \gamma}$.
Then $t^{-e_i} a_{ij} \in \OO_{X, \gamma}$ for all $j$.
Thus $t^{-e_i} \omega'_i \in E_{\gamma}$ and
$t^{-e_i} \omega'_i \not= 0 $ in $E_{\gamma} \otimes \kappa(\gamma)$.
Let $\Gamma$ be the Zariski closure of $\{ \gamma \}$.
If we choose a general closed point $x_0$ of $\Gamma$, then
$\omega'_i \not= 0$  in $E'_{x_0} \otimes \kappa(x_0)$ and
$t^{-e_i} \omega'_i \not= 0$ in $E_{x_0} \otimes \kappa(x_0)$.
On the other hand, there is an open neighborhood $U_{x_0}$ of $x_0$ such that
\[
h(t^{-e_i} \omega'_i, t^{-e_i} \omega'_i)(x) = h'(t^{-e_i} \omega'_i, t^{-e_i} \omega'_i)(x)
\]
for $x \in U_{x_0} \cap U$.
Thus if we set
\[
f(x) = h(t^{-e_i} \omega'_i, t^{-e_i} \omega'_i)(x) = \vert t \vert^{-2e_i}
h'(\omega'_i, \omega'_i)(x)
\]
on $U_{x_0} \cap U$,
then $\lim_{x \to x_0} f(x) = h(t^{-e_i} \omega'_i, t^{-e_i} \omega'_i)(x_0) = 0$
because $t = 0$ at $x_0$.
This is a contradiction because $t^{-e_i} \omega'_i \not= 0$ in $E_{x_0} \otimes \kappa(y)$.
Therefore we can see that $a_{ij} \in \OO_{X,\gamma}$ for all $i, j$.
In the same way, $b_{ij} \in \OO_{X,\gamma}$ for all $i, j$.

\medskip
By the above claim, $\{ \omega_1, \ldots, \omega_r \}$ and
$\{ \omega'_1, \ldots, \omega'_r \}$ generate the same $\OO_{X, \gamma}$-module
in $V$. Thus $E_{\gamma} = E'_{\gamma}$. Hence we get $E = E'$.

Let $x$ be an arbitrary closed point of $X$.
Let $v, v' \in E_x \otimes \kappa(x)$.
Choose $\omega, \omega' \in E_x$ such that
$\omega$ and $\omega'$ give rise to $v$ and $v'$ in  $E_x \otimes \kappa(x)$.
Then there is a neighborhood $U_x$ of $x$ such that
$h(\omega, \omega')(y) = h'(\omega, \omega')(y)$ for all $y \in U_x \cap U$.
Thus
\[
h(\omega, \omega')(x) = \lim_{y \to x} h(\omega, \omega')(y)
= \lim_{y \to x} h'(\omega, \omega')(y) = h'(\omega, \omega')(x),
\]
which means that $h_x(v, v') = h'_x(v, v')$.
\QED

\begin{Proposition}
\label{prop:loc:int:det:H}
Let $(E, h)$ be a $C^{\infty}$-hermitian locally free coherent
sheaf on $X$. 
Let $x_1, \ldots, x_r$ be a $K$-linearly independent elements of $E_{\eta}$.
Then $\log(\det (h(x_i, x_j)))$ is a locally integrable function.
\end{Proposition}

\Proof
Let $W$ be a vector subspace of $E_{\eta}$ generated by $x_1, \ldots, x_r$.
By Proposition~\ref{prop:subsheaves:vector:subspaces}, there is 
a saturated $\OO_X$-subsheaf $F$ of $E$ with $F_{\eta} = W$.
First we assume that $F$ and $E/F$ are locally free.
For a closed point $x \in X$, let $\{ \omega_1, \ldots, \omega_r \}$
be a local basis of $F_x$. Then we can find a matrix $A = (a_{ij})$ such that
$a_{ij} \in K$ for all $i, j$ and
$x_i = \sum_{j=1}^r a_{ij}\omega_j$ for all $i$.
Then
\[
\det (h(x_i, x_j)) = \vert \det(A) \vert^2 \det (h(\omega_i, \omega_j)).
\]
Since $F$ and $E/F$ are locally free, $\det (h(\omega_i, \omega_j))$
is a non-zero $C^{\infty}$-function around $x$ and $\det(A)$ is a non-zero
rational function on $X$. Thus $\log(\det (h(x_i, x_j)))$ is locally integrable around $x$.

\medskip
In general, if we set $Q = E/F$, then there is a proper birational morphism
$\mu : Y \to X$ of smooth algebraic varieties over $\CC$ such that
\[
\mu^*(Q)/(\text{the torsion part of $\mu^*(Q)$})
\]
is locally free. 
We set $F'  = \Ker(\mu^*(E) \to \mu^*(Q)/(\text{the torsion part of $\mu^*(Q)$}))$.
Then $F'$ and $\mu^*(E)/F'$ are locally free.
Thus, since $F'_{\eta} = W$,
\[
\log (\det (\mu^*(h)(x_i, x_j))) = \mu^*(\log (\det (h(x_i, x_j))))
\]
is a locally integrable function on $Y$.
Therefore so is $\log \det (h(x_i, x_j))$ on $X$ by virtue of
\cite[Proposition~1.2.5]{LMSemi}
\QED

\subsection{Arakelov geometry}
\setcounter{Theorem}{0}
For basic definitions concerning Arakelov geometry, we refer to \cite[Section~1]{MoArht}.
Let $X$ be a projective arithmetic variety.
We use several kinds of positivity of a $C^{\infty}$-hermitian invertible sheaf on $X$ (like ampleness, nefness and bigness)
as defined in \cite[Section~2]{MoArht}.
Let $\overline{H} = (\overline{H}_1, \ldots, \overline{H}_d)$ be a sequence of 
nef $C^{\infty}$-hermitian invertible sheaves on $X$, where $d = \dim X_{\QQ}$. Note that
the sequence is empty in the case of $d = 0$. We say $\overline{H}$ is {\em fine} if
$(X; \overline{H}_1, \ldots, \overline{H}_d)$ gives rise to a fine polarization of the function field of $X$
(for details, see \cite[Section~6.1]{MoCycle}). For example, if $\overline{H}_i$'s are nef and big,
then $\overline{H}$ is fine.
Finally we consider the following lemma.

\begin{Lemma}
\label{lem:isom:chow:codim:two}
Let $X$ be a generically smooth arithmetic variety and $U$ a Zariski open set of $X$ with
$\codim(X \setminus U) \geq 2$. Then the natural homomorphism
\[
\aChow^1_D(X) \to \aChow^1_D(U)
\]
is injective.
\end{Lemma}

\Proof
Let $(D, T)$ be an arithmetic cycle of codimension one on $X$.
We assume that $(\rest{D}{U}, \rest{T}{U}) = \widehat{(\rest{\phi}{U})}$
for some non-zero rational function $\phi$ on $X$. Then,
since $\codim(X \setminus U) \geq 2$, we have
$(D, T) = \widehat{(\phi)}$.
\QED

\renewcommand{\theTheorem}{\arabic{section}.\arabic{Theorem}}
\renewcommand{\theClaim}{\arabic{section}.\arabic{Theorem}.\arabic{Claim}}
\renewcommand{\theequation}{\arabic{section}.\arabic{Theorem}.\arabic{Claim}}

\section{Birationally $C^{\infty}$-hermitian torsion free coherent sheaves \\
on a normal arithmetic variety}

Let $X$ be a normal arithmetic variety.
Let $E$ be a torsion free coherent sheaf on $X$.
We say a pair $(E, h)$ is called a 
{\em birationally $C^{\infty}$-hermitian torsion free coherent sheaf} on $X$
if there are a proper birational morphism $\mu : X' \to X$ of 
normal arithmetic varieties,
a $C^{\infty}$-hermitian locally free coherent sheaf $(E', h')$ on $X'$,
and a Zariski open set $U$ of $X$
with the following properties: 
\begin{enumerate}
\renewcommand{\labelenumi}{(\arabic{enumi})}
\item
$X'$ and $U$ are generically smooth.

\item
$\codim(X \setminus U) \geq 2$.

\item
$\mu : X' \to X$ is an isomorphism over $U$, that is,
if we set $U' = \mu^{-1}(U)$, then
$\rest{\mu}{U'} : U' \overset{\sim}{\longrightarrow} U$.

\item
$E$ is locally free on $U$ and
$h$ is a $C^{\infty}$-hermitian metric of $\rest{E}{U}$ over $U(\CC)$.

\item
$(\rest{\mu}{U'})^*(\rest{(E, h)}{U})$ is isometric to
$\rest{(E', h')}{U'}$.
\end{enumerate}
This $C^{\infty}$-hermitian locally free coherent sheaf $(E', h')$
is called a {\em model of $(E, h)$ in terms of $\mu : X' \to X$}.
Note that if $\mu' : X'' \to X'$ is a proper birational morphism
of normal and generically smooth arithmetic varieties, then
${\mu'}^*(E', h')$ is also a model of $(E, h)$ in terms of $\mu \circ \mu' : X'' \to X$.
For, let $X'_0$ be the maximal Zariski open set over which
$\mu'$ is an isomorphism. Then $\codim(X' \setminus X_0)\geq 2$.
Thus if we set $V = \mu(U' \cap X'_0)$, then we can see the above properties for $V$.

\begin{Proposition}
\label{prop:saturated:subsheaf:birat:herm}
Let $X$ be a normal arithmetic variety and
$(E, h)$ a birationally
$C^{\infty}$-hermitian torsion free coherent sheaf on $X$.
Let $F$ be a saturated $\OO_X$-subsheaf of $E$.
Let $h_{F\hookrightarrow E}$ \rom{(}resp. $h_{E \twoheadrightarrow E/F}$\rom{)}
be the submetric of $F$ induced by
$F \hookrightarrow E$ and $h$ 
\rom{(}resp. the quotient metric of $E/F$ induced by
$E \twoheadrightarrow E/F$ and $h$\rom{)}
on a big Zariski open set of $X$, i.e.,
a Zariski open set whose complement has the codimension greater than or equal to
$2$.
Then $(F, h_{F\hookrightarrow E})$ and $(E/F, h_{E \twoheadrightarrow E/F})$
are also a birationally
$C^{\infty}$-hermitian torsion free coherent sheaf on $X$.
\end{Proposition}

\Proof
Let  $\eta$ be the generic point of $X$.
Let $(E', h')$ be a model of $(E, h)$ in terms of
$\mu : X' \to X$.
Let $F'$ be a saturated $\OO_{X'}$-subsheaf $F'$ of $E'$ with
$F'_{\eta} = F_{\eta}$ \rom{(}cf.  Proposition~\rom{\ref{prop:subsheaves:vector:subspaces}}\rom{)}.
We set $Q = E'/F'$.
By \cite[Theorem~1 in Chapter~4]{RayFlat},
there is a proper birational morphism $\mu' : X'' \to X'$
of normal and generically smooth arithmetic varieties
such that ${\mu'}^*(Q)/(\text{torsion})$ is locally free.
Let 
\[
F'' = \Ker ({\mu'}^*(E') \to {\mu'}^*(Q)/(\text{torsion})).
\]
Then $F''$ and ${\mu'}^*(E')/F''$ are locally free.
Thus
\[
(F'', {\mu'}^*(h')_{F'' \hookrightarrow {\mu'}^*(E')})
\quad\text{and}\quad
({\mu'}^*(E')/F'', {\mu'}^*(h')_{{\mu'}^*(E') \twoheadrightarrow {\mu'}^*(E')/F''})
\]
yield models
of $(F, h_{F\hookrightarrow E})$ and $(E/F, h_{E \twoheadrightarrow E/F})$
respectively
because ${\mu'}^*(E', h')$ gives rise to  a model of $(E, h)$.
\QED

\begin{Proposition}
\label{prop:well:def:deg:H}
We assume that $X$ is projective.
Let $\overline{H} = (\overline{H}_1, \ldots, \overline{H}_d)$ be
a sequence of nef $C^{\infty}$-hermitian invertible sheaves on $X$,
where $d = \dim X_{\QQ}$.
Then the quantity 
\[
\adeg(\acherncl_1(\mu^*(\overline{H}_1)) \cdots
\acherncl_1(\mu^*(\overline{H}_d)) \cdot \acherncl_1(E', h') )
\]
does not depend on the choice of a model
$(E', h')$ in terms of $\mu : X' \to X$. It is denoted by
$\adeg_{\overline{H}}(E, h)$ and is called the arithmetic degree of
$(E, h)$ with respect to $\overline{H}$.
\end{Proposition}

\Proof
Let us begin with the following lemma.

\begin{Lemma}
\label{lem:equal:two:arith:degree}
Let $\nu : Y \to X$ be a birational morphism of normal and projective arithmetic varieties
such that $Y$ is generically smooth.
Let $(E, h)$ and $(E', h')$ be $C^{\infty}$-hermitian locally free coherent sheaves on $Y$.
We assume that there is a Zariski open set $U$ of $X$ such that
$\codim(X \setminus U) \geq 2$ and
$\nu$ is an isomorphism over $U$, that is,
if we set $V = \nu^{-1}(U)$, then $\rest{\nu}{V} : V \overset{\sim}{\longrightarrow} U$.
Let $\overline{L}_1, \ldots, \overline{L}_d$ be $C^{\infty}$-hermitian invertible
sheaves on $X$, where $d = \dim X_{\QQ}$.
If $\rest{(E, h)}{V}$ is isometric to $\rest{(E', h')}{V}$, then
\begin{multline*}
\adeg( \acherncl_1(\nu^*(\overline{L}_1)) \cdots
\acherncl_1(\nu^*(\overline{L}_d)) \cdot \acherncl_1(E, h) ) \\
=
\adeg( \acherncl_1(\nu^*(\overline{L}_1)) \cdots
\acherncl_1(\nu^*(\overline{L}_d)) \cdot \acherncl_1(E', h')).
\end{multline*}
\end{Lemma}

\Proof
Let $\eta$ be the generic point of $Y$ and $x_1, \ldots, x_r$ a basis of
$E_{\eta}$. Let $x'_1, \ldots, x'_r$ be the corresponding basis of $E'_{\eta}$
with $x_1, \ldots, x_r$.
Let $Y^{(1)}$ be the set of all codimension one points of $Y$.
Then $\acherncl_1(E, h)$ and $\acherncl_1(E', h')$ are represented by
\[
\left( \sum_{\gamma \in Y^{(1)} }
\ell_{\OO_{Y, \gamma}}(E; x_1, \ldots, x_r)\overline{\{ \gamma \}},
\ - \log (\det (h(x_i, x_j))) \right)
\]
and
\[
\left( \sum_{\gamma \in Y^{(1)} }
\ell_{\OO_{Y, \gamma}}(E'; x'_1, \ldots, x'_r)\overline{\{ \gamma \}},
\ - \log (\det (h'(x'_i, x'_j))) \right)
\]
respectively.
By Proposition~\ref{prop:coincide:hermitian:locally:free:sheaf},
we can see that 
\[
\det (h(x_i, x_j)) = \det (h'(x'_i, x'_j))
\]
on $Y(\CC)$.
Here
\[
\ell_{\OO_{Y, \gamma}}(E; x_1, \ldots, x_r) = 
\ell_{\OO_{Y, \gamma}}(E'; x'_1, \ldots, x'_r)
\]
for all $\gamma \in V^{(1)}$.
Moreover,  for  $\gamma \in Y^{(1)} \setminus V^{(1)}$,
since $\codim(\nu(\overline{\{ \gamma\} })) \geq 2$,
\[
\adeg ( \acherncl_1(\nu^*(\overline{L}_1)) \cdots
\acherncl_1(\nu^*(\overline{L}_d)) \cdot (\overline{\{ \gamma\} }, 0)) = 0
\]
by the projection formula (cf. \cite[Proposition~1.2 and Proposition~1.3]{MoArht}).
Thus we have our lemma.
\QED

\bigskip
Let us go back to the proof of Proposition~\ref{prop:well:def:deg:H}.
Let $(E_1, h_1)$ and $(E_2, h_2)$
be two models of $(E, h)$ in terms of
$\mu_1 : X_1 \to X$ and $\mu_2 : X_2 \to X$ respectively.
We can choose a normal, projective and generically smooth arithmetic variety $Y$ and
birational morphisms $\pi_1 : Y \to X_1$ and $\pi_2 : Y \to X_2$ with
$\mu_1 \circ \pi_1 = \mu_2 \circ \pi_2$.
We set $\nu = \mu_1 \circ \pi_1 = \mu_2 \circ \pi_2$.
First of all, by  the projection formula, we have
\begin{multline*}
\adeg( \acherncl_1(\mu_1^*(\overline{H}_1)) \cdots
\acherncl_1(\mu_1^*(\overline{H}_d)) \cdot \acherncl_1(E_1, h_1)) \\
=
\adeg(\acherncl_1(\nu^*(\overline{H}_1)) \cdots
\acherncl_1(\nu^*(\overline{H}_d)) \cdot \acherncl_1(\pi_1^*(E_1, h_1)))
\end{multline*}
and
\begin{multline*}
\adeg( \acherncl_1({\mu_2}^*(\overline{H}_1)) \cdots
\acherncl_1({\mu_2}^*(\overline{H}_d)) \cdot \acherncl_1(E_2, h_2))\\
 =
\adeg( \acherncl_1(\nu^*(\overline{H}_1)) \cdots
\acherncl_1(\nu^*(\overline{H}_d)) \cdot \acherncl_1(\pi_2^*(E_2, h_2))).
\end{multline*}
Moreover, by Lemma~\ref{lem:equal:two:arith:degree},
\begin{multline*}
\adeg( \acherncl_1(\nu^*(\overline{H}_1)) \cdots
\acherncl_1(\nu^*(\overline{H}_d)) \cdot \acherncl_1(\pi_1^*(E_1, h_1)) ) \\
=
\adeg( \acherncl_1(\nu^*(\overline{H}_1)) \cdots
\acherncl_1(\nu^*(\overline{H}_d)) \cdot \acherncl_1(\pi_2^*(E_2, h_2))).
\end{multline*}
Thus we get the assertion.
\QED

Let $X$ be a normal arithmetic variety and $(E, h)$ a birationally $C^{\infty}$-hermitian
torsion free sheaf on $X$.
Let $\pi : X' \to X$ be a proper birational morphism of
normal arithmetic varieties and $(E', h')$ a birationally $C^{\infty}$-hermitian
torsion free sheaf on $X'$.
We say $(E, h)$ is {\em birationally dominated by
$(E', h')$ by means of $\pi : X' \to X$} if there is a Zariski open set $U$ of $X$ with
the following properties:
\begin{enumerate}
\renewcommand{\labelenumi}{(\arabic{enumi})}
\item
$\codim(X \setminus U) \geq 2$ and $U$ is generically smooth.

\item
$(E, h)$ is a $C^{\infty}$-hermitian locally free sheaf over $U$.

\item
If we set $U' = \pi^{-1}(U)$, then $\rest{\pi}{U'} : U' \overset{\sim}{\longrightarrow} U$.

\item
$(\rest{\pi}{U'})^*(\rest{(E, h)}{U})$ is isometric to
$\rest{(E', h')}{U'}$.
\end{enumerate}

Then we have the following:

\begin{Proposition}
\label{prop:birat:dominant}
The notation is the same as above. We assume that
$(E, h)$ is birationally dominated by
$(E', h')$ by means of $\pi : X' \to X$.
\begin{enumerate}
\renewcommand{\labelenumi}{(\arabic{enumi})}
\item
Let $F$ be a saturated $\OO_X$-subsheaf of $E$
and $F'$ the corresponding saturated $\OO_{X'}$-subsheaf of $E'$
with $F$.
Then $(F, h_{F \hookrightarrow E})$ and
$(E/F, h_{E \twoheadrightarrow E/F})$ are birationally dominated by
$(F', h'_{F' \hookrightarrow E'})$ and
$(E'/F', h'_{E' \twoheadrightarrow E'/F'})$ respectively.

\item
We assume that $X$ and $X'$ are projective.
Let $\overline{H} = (\overline{H}_1, \ldots, \overline{H}_d)$ be
a sequence of nef $C^{\infty}$-hermitian invertible sheaves on $X$,
where $d = \dim X_{\QQ}$.
Then $\adeg_{\overline{H}}(E, h) = \adeg_{\pi^*(\overline{H})}(E', h')$.
\end{enumerate}
\end{Proposition}

\Proof
(1) There is a Zariski open set $U_1$ such that
$U_1 \subseteq U$, $\codim(X \setminus U_1) \geq 2$ and that
$\rest{E}{U_1}$ and $\rest{E/F}{U_1}$ are locally free.
We set $U'_1 = \pi^{-1}(U_1)$. Then
$(\rest{\pi}{U'})^*(\rest{(F, h_{F \hookrightarrow E})}{U_1})$ is isometric to
$\rest{(F', h'_{F' \hookrightarrow E'})}{U'_1}$.
Thus our assertions follow.

\medskip
(2) Let $(E'', h'')$ be a model of $(E', h')$ in terms of a birational morphism
$\mu : Y \to X'$. Then it is easy to see that $(E'', h'')$ is a model of $(E, h)$
in terms of $\pi \circ \mu : Y \to X$.
Thus we have (2) by Proposition~\ref{prop:well:def:deg:H}.
\QED

\section{Finiteness of subsheaves with bounded arithmetic degree}

In this section, we would like to give the proof of the main theorem of this note.

\begin{Theorem}
\label{thm:finite:subsheaf}
Let $X$ be a normal projective arithmetic variety and
$(E, h)$ a birationally $C^{\infty}$-hermitian torsion free coherent sheaf on $X$.
Let $\overline{H} = (\overline{H}_1, \ldots, \overline{H}_d)$ be
a fine sequence of nef $C^{\infty}$-hermitian invertible sheaves on $X$,
where $d = \dim X_{\QQ}$.
For any real number $c$, the set of all non-zero saturated $\OO_{X}$-subsheaf $F$ of $E$
with $\adeg_{\overline{H}}(\acherncl_1(F, h_{F\hookrightarrow E})) \geq c$ is finite, where
$h_{F\hookrightarrow E}$ is the submetric of $F$ induced by $h$ over a big open set.
\end{Theorem}

\Proof
Let $(E', h')$ be a model of $(E, h)$ in terms of $\mu : X' \to X$.
Let $\eta$ be the generic point of $X$.
For each vector subspace $W$ of $E_{\eta}$,
let $F$ (resp. $F'$) be a saturated $\OO_X$-subsheaf of $E$
(resp. $\OO_{X'}$-subsheaf of $E'$) induced by $W$.
Then, by Proposition~\ref{prop:birat:dominant},
\[
\adeg_{\overline{H}}(F, h_{F \hookrightarrow E}) =
\adeg_{\mu^*(\overline{H})}(F', h_{F' \hookrightarrow E'}).
\]
Therefore we may assume that $X$ is generically smooth,
$E$ is locally free and $h$ is a $C^{\infty}$-hermitian metric of $E$.

For each $0 < s < \rank E$, let
$\Sigma_s(X, E)$ be the set of 
all saturated rank $s$ $\OO_X$-subsheaves of $E$.
First let us see that, for any real number $c$, the set
\[
\{ L \in \Sigma_1(X, E) \mid 
\adeg_{\overline{H}}(F, h_{F\hookrightarrow E}) \geq c \}
\]
is finite.
Let $\pi : P = \Proj(\bigoplus_{d \geq 0} \Sym^d(E^{\vee})) \to X$
be the projective bundle and $\OO_{P}(1)$ the tautological line
bundle of $P$.
Let $h_P$ be the quotient hermitian metric of $\OO_P(1)$ by
using the surjective homomorphism
$\pi^*(E^{\vee}) \to \OO_{P}(1)$ and the hermitian metric $\pi^*(h^{\vee})$.
In other words,
the metric $h_P^{-1}$ of $\OO_P(-1)$ is the submetric
induced by the injective homomorphism $\OO_P(-1) \to \pi^*(E)$ and $\pi^*(h)$
(cf. (3) of Proposition~\ref{prop:dual:hermitian:vector}).
Let $P_{\eta}$ be the generic fiber of $\pi : P \to X$, and
$K$ the function field of $X$.

For a $K$-rational point $x$ of $P_{\eta}$,
let us introduce $\Delta_x$, $U_x$, $V_x$ and $s_x$ as follows:
$\Delta_x$ is the Zariski closure of $x$ in $P$ and
$U_x$ is the maximal open set of $X$ over which
$\rest{\pi}{\Delta_x} : \Delta_x \to X$ is an isomorphism.
Further $V_x = (\rest{\pi}{\Delta_x})^{-1}(U_x)$ and
$s_x : U_x \to P$ is the section induced by  the isomorphism
$\rest{\pi}{V_x} : V_x \to U_x$

Let $\Sigma_1(K, E_{\eta})$ be the set of all $1$-dimensional vector subspaces of $E_{\eta}$
over $K$.
Then, by Proposition~\ref{prop:invertible:section:1to1},
there is a natural bijection
\[
P_{\eta}(K) \to \Sigma_1(K, E_{\eta}).
\]
Moreover let $\Sigma_1(X, E)$ be the set of 
all saturated rank one $\OO_X$-subsheaves of $E$.
By Proposition~\ref{prop:subsheaves:vector:subspaces},
we have a bijective map
\[
\Sigma_1(X, E) \to \Sigma_1(K, E_{\eta}).
\]
Therefore there is a natural bijection between
$P_{\eta}(K)$ and $\Sigma_1(X, E)$.
For a $K$-rational point $x$ of $P_{\eta}$,
the corresponding 
saturated rank one $\OO_X$-subsheaf of $E$ is denoted by  $L(x)$.
Then, by using Proposition~\ref{prop:invertible:section:1to1},
we can see that
$L(x)$ has the following property:
Let $s_x^*(\OO_P(-1)) \to s_x^*\pi^*(E) = \rest{E}{U_x}$ be
the homomorphism from the natural homomorphism
$\OO_P(-1) \to \pi^*(E)$ by applying $s_x^*$.
Then the image of $s_x^*(\OO_P(-1)) \to \rest{E}{U_x}$ is
$\rest{L(x)}{U_x}$.
Let $h_x$ be the submetric of $L(x)$ induced by $h$.

\begin{Claim}
$\acherncl_1(L(x), h_x) = 
(\rest{\pi}{\Delta_x})_* \left(\acherncl_1\left(\rest{(\OO_P(-1), h_P^{-1})}{\Delta_x}\right)\right)$.
\end{Claim}

Since the metric $h_P^{-1}$ is the submetric of  $\OO_P(-1)$ induced
by $\pi^*(h)$, we can see that $s_x^*(\OO_P(-1), h_P^{-1})$ is isometric
to $\rest{(L(x), h_x)}{U_x}$. Thus
$\rest{(\OO_P(-1), h_P^{-1})}{V_x}$ is
isometric to $(\rest{\pi}{V_x})^*(\rest{(L(x), h_x)}{U_x})$,
which implies that
\begin{align*}
(\rest{\pi}{V_x})_* \left(\acherncl_1\left(\rest{(\OO_P(-1), h_P^{-1})}{V_x}\right)\right) & =
(\rest{\pi}{V_x})_* \left(\acherncl_1\left(  (\rest{\pi}{V_x})^*(\rest{(L(x), h_x)}{U_x}) \right)\right) \\
&=\acherncl_1(\rest{(L(x), h_x)}{U_x}).
\end{align*}
This means that the assertion of the claim holds over $U_x$. Thus so does over $X$
by Lemma~\ref{lem:isom:chow:codim:two}.

\medskip
For a $K$-rational point $x$ of $P_{\eta}$,  the height $h_{\OO(1)}(x)$ with respect to
$\OO_P(1)$ and $(X, \overline{H})$ is given by
\[
h_{\OO(1)}(x) = \adeg\left( \acherncl_1((\rest{\pi}{\Delta_x})^*(\overline{H}_1)) 
\cdots  \acherncl_1((\rest{\pi}{\Delta_x})^*(\overline{H}_d))  \cdot
\acherncl_1\left(\rest{(\OO_P(1), h_P)}{\Delta_x}\right) \right).
\]
By using the above claim and the projection formula,
\begin{align*}
-h_{\OO_P(1)}(x) & =
\adeg\left( \acherncl_1((\rest{\pi}{\Delta_x})^*(\overline{H}_1)) \cdots
\acherncl_1((\rest{\pi}{\Delta_x})^*(\overline{H}_d)) \cdot
\acherncl_1\left(\rest{(\OO_P(-1), h_P^{-1})}{\Delta_x}\right) \right) \\
& = \adeg\left( \acherncl_1(\overline{H_1}) \cdots
\acherncl_1(\overline{H}_d) \cdot \acherncl_1(L(x), h_x)\right) =
\adeg_{\overline{H}}(L(x), h_x).
\end{align*}
Thus we have a bijective correspondence between
\[
\{ L \in \Sigma_1(X, E) \mid 
\adeg_{\overline{H}}(F, h_{F\hookrightarrow E}) \geq c \}
\]
and
\[
\{ x \in P_{\eta}(K) \mid h(x) \leq -c \}.
\]
On the other hand, by virtue of Northcott's theorem
over finitely generated field (cf. \cite[Theorem~4.3]{MoArht}),
$\{ x \in P_{\eta}(K) \mid h(x) \leq -c \}$ is a finite set.
Therefore we get the case where $s = 1$.

\medskip
For $F \in \Sigma_s(X, E)$, let $\lambda(F)$ be the saturation of
\[
\bigwedge^s F/(\text{the torsion part of $\bigwedge^s F$})
\]
in $\bigwedge^s E$.

\begin{Claim}
If $\lambda(F) = \lambda(F')$, then $F = F'$.
\end{Claim}

We assume that $\lambda(F) = \lambda(F')$.
Let $K$ be the function field of $X$.
Then, using Pl\"ucker coordinates over $K$,
we can see that $F \otimes K = F' \otimes K$.
Thus, by Lemma~\ref{unique:saturate:subsheaf}, $F' = F$.

\medskip
Let $h_{\lambda(F)} = (\bigwedge^s h)_{\lambda(F) \hookrightarrow \bigwedge^s E}$. 
Then, by Proposition~\ref{prop:hermitian:metric:wedge},
\[
\acherncl_1(F, h_F) = \acherncl_1(\lambda(F), h_{\lambda(F)}).
\]
Therefore,
by using the above claim and
the case where $s=1$, our theorem follows.
\QED

Let $X$ be a normal and projective arithmetic variety and
$(E, h)$ a birationally $C^{\infty}$-hermitian torsion free coherent sheaf on $X$.
Let $\overline{H} = (\overline{H}_1, \ldots, \overline{H}_d)$ be
a fine sequence of nef $C^{\infty}$-hermitian invertible sheaves on $X$.
For a non-zero saturated $\OO_X$-subsheaf $G$ of $E$, we set
\[
\hat{\mu}_{\overline{H}}(G, h_{G \hookrightarrow E}) = 
\frac{\adeg_{\overline{H}}(G, h_{G \hookrightarrow E})}{\rank G}.
\]
A saturated $\OO_X$-subsheaf $F$ of $E$ is called a {\em maximal slope sheaf of $(E, h)$
with respect to $\overline{H}$} if
$\hat{\mu}_{\overline{H}}(F, h_{F \hookrightarrow E})$ gives rise to the maximal value of the set
\[
\left\{ \hat{\mu}_{\overline{H}}(G, h_{G \hookrightarrow E})
\mid \text{$G$ is a non-zero saturated $\OO_X$-subsheaf of $E$}\right\}.
\]
Moreover a maximal slope sheaf $F$ of $(E, h)$ is called
a {\em maximal destabilizing sheaf of $(E, h)$  with respect to $\overline{H}$} if
$\rank F$ is maximal among all maximal slope sheaves of $(E, h)$.
As a corollary of Theorem~\ref{thm:finite:subsheaf}, we have the following:

\begin{Corollary}
\label{cor:max:destabilizing:sheaf}
There is a maximal destabilizing sheaf of $(E, h)$ with respect to $\overline{H}$.
\end{Corollary}

\section{Arithmetic first Chern class of a subsheaf}
Let $X$ be a normal and generically smooth 
arithmetic variety and $\eta$ the generic point of $X$.
Let $(E, h)$ be a $C^{\infty}$-hermitian locally free sheaf on $X$.
Let $F$ be an $\OO_X$-subsheaf of $E$.
Let $x_1, \ldots, x_r$ be a basis of $F_{\eta}$.
Let us consider an arithmetic codimension one cycle $z(F; x_1, \ldots, x_r)$
(i.e., an element of $ \in \widehat{Z}^1_D(X)$) given by
\[
z(F; x_1, \ldots, x_r) = \left( \sum_{\Gamma} 
\ell_{\OO_{X, \Gamma}}(F_{\Gamma}; x_1, \ldots, x_r)\Gamma, 
-\log \det (h(x_i, x_j)) \right).
\]
Note that $\log \det (h(x_i, x_j))$ is locally integrable on $X(\CC)$ by
Proposition~\ref{prop:loc:int:det:H}.
Let $x'_1, \ldots, x'_r$ be another  basis of $F_{\eta}$.
There is an $r \times r$-matrix $A = (a_{ij})$ with
$x'_i = \sum_{j=1}^r a_{ij} x_j$.
Using (2) of Corollary~\ref{cor:base:change:rule},
we can see that
\[
z(F; x'_1, \ldots, x'_r) = z(F; x_1, \ldots, x_r) + \widehat{(\det(A))}.
\]
Therefore the class of $z(F; x_1, \ldots, x_r)$ in $\aChow^1_D(X)$ does not
depend on the choice of $x_1, \ldots, x_r$.
We denote the class of $z(F; x_1, \ldots, x_r)$ in $\aChow^1_D(X)$
by $\acherncl_1(F \hookrightarrow E, h)$.
If $F = E$, then $\acherncl_1(E \hookrightarrow E, h)$ is equal to the usual
$\acherncl_1(E, h)$.
Note that
\[
\acherncl_1(F \hookrightarrow E, h) = 
\acherncl_1(F, h_{F\hookrightarrow E})
\]
if $F$ is saturated in $E$.
More generally, we have the following:

\begin{Proposition}
\label{prop:diff:saturation}
Let $F$ be an $\OO_X$-subsheaf of $E$ and $\widetilde{F}$ the saturation of $F$ in $E$.
Then $\acherncl_1(\widetilde{F}, h_{\widetilde{F}\hookrightarrow E}) - 
\acherncl_1(F \hookrightarrow E, h)$
is represented by an arithmetic divisor
\[
\left( \sum_{\text{\rom{$\Gamma$ : prime divisor}}} \length_{\OO_{X, \Gamma}}(\widetilde{F}_{\Gamma}/F_{\Gamma}) \Gamma, 0 \right).
\]
In particular, if $\overline{H} = (\overline{H}_1, \ldots, \overline{H}_d)$
is a sequence of nef $C^{\infty}$-hermitian invertible sheaves on $X$,
then
\[
\adeg(\acherncl_1(\overline{H}_1) \cdots \acherncl_1(\overline{H}_d) \cdot \acherncl_1(F \hookrightarrow E, h)) \leq
\adeg(\acherncl_1(\overline{H}_1) \cdots \acherncl_1(\overline{H}_d) \cdot  \acherncl_1(\widetilde{F}, h_{\widetilde{F}\hookrightarrow E})).
\]
\end{Proposition}

\Proof
Let $\eta$ be the generic point of $X$.
Let $\{ x_1, \ldots, x_r \}$ be a basis of $F_{\eta}$.
Then $\{ x_1, \ldots, x_r\}$ also gives rise to a basis of $\widetilde{F}_{\eta}$.
Thus $\acherncl_1(\widetilde{F}, h_{\widetilde{F}\hookrightarrow E}) - \acherncl_1(F \hookrightarrow E, h)$
is represented by
\[
\left( \sum_{\Gamma}(\ell_{\OO_{X, \Gamma}}(\widetilde{F}_{\Gamma}; x_1, \ldots, x_r) -
\ell_{\OO_{X, \Gamma}}(F_{\Gamma}; x_1, \ldots, x_r) ) \Gamma, 0 \right).
\]
Hence it is sufficient to see that
\[
\ell_{\OO_{X, \Gamma}}(\widetilde{F}_{\Gamma}; x_1, \ldots, x_r) -
\ell_{\OO_{X, \Gamma}}(F_{\Gamma}; x_1, \ldots, x_r) = 
\length_{\OO_{X, \Gamma}}(\widetilde{F}_{\Gamma}/F_{\Gamma}) 
\]
for all $\Gamma$.
Let $a$ be an element of $\OO_{X,\Gamma} \setminus \{ 0 \}$ such that
$a x_i \in \OO_{X,\Gamma}$ for all $i$.
Then
\begin{align*}
\ell_{\OO_{X, \Gamma}}(\widetilde{F}_{\Gamma}; x_1, \ldots, x_r) & = \length_{\OO_{X,\Gamma}}(
\widetilde{F}_{\Gamma}/\OO_{X,\Gamma} ax_1 + \cdots + \OO_{X,\Gamma} ax_r) - r \ord_{\Gamma}(a), \\
\ell_{\OO_{X, \Gamma}}(F_{\Gamma}; x_1, \ldots, x_r) & = \length_{\OO_{X,\Gamma}}(
F_{\Gamma}/\OO_{X,\Gamma} ax_1 + \cdots + \OO_{X,\Gamma} ax_r) - r \ord_{\Gamma}(a).
\end{align*}
Therefore we get our proposition.
\QED

\section{Arithmetic Harder-Narasimham filtration}
Let $X$ be a normal and projective arithmetic variety
and $\overline{H} = (\overline{H}_1, \ldots, \overline{H}_d)$ a fine
sequence of nef $C^{\infty}$-hermitian invertible sheaves.
Let $(E, h)$ be a birationally $C^{\infty}$-hermitian torsion free
coherent sheaf on $X$. 
$(E, h)$ is said to be {\em arithmetically $\mu$-semistable} with respect to $\overline{H}$ if,
for any non-zero saturated $\OO_{X}$-subsheaf $F$ of $E$,
\[
\hat{\mu}_{\overline{H}}( F, h_{F\hookrightarrow E})
\leq 
\hat{\mu}_{\overline{H}}(E, h).
\]
A filtration
\[
0 = E_0 \subsetneq E_1 \subsetneq \cdots \subsetneq E_l = E
\]
of $\OO_X$-subsheaves of $E$ is called a {\em saturated filtration of $E$}
if $E_i/E_{i-1}$ is torsion free for every $1 \leq i \leq l$.
Moreover we say a saturated filtration
$0 = E_0 \subsetneq E_1 \subsetneq \cdots \subsetneq E_l = E$
of $E$ is an {\em arithmetic Harder-Narasimham filtration of $(E, h)$
with respect to $\overline{H}$} if
\begin{enumerate}
\renewcommand{\labelenumi}{(\arabic{enumi})}
\item
Let $h_{E_i/E_{i-1}}$ be a $C^{\infty}$-hermitian metric of $E_i/E_{i-1}$
induced by $h$, that is,
\[
h_{E_i/E_{i-1}} = (h_{E_i \hookrightarrow E})_{E_i \twoheadrightarrow E_i/E_{i-1}} = 
(h_{E \twoheadrightarrow E/E_{i-1}})_{E_i/E_{i-1}  \hookrightarrow E/E_{i-1}}.
\]
Then $(E_i/E_{i-1}, h_{E_i/E_{i-1}})$ is arithmetically $\mu$-semistable with respect to $\overline{H}$.

\item
$\hat{\mu}_{\overline{H}}(E_1/E_0, h_{E_1/E_0}) >
\hat{\mu}_{\overline{H}}(E_2/E_1, h_{E_2/E_1}) > \cdots >
\hat{\mu}_{\overline{H}}(E_l/E_{l-1}, h_{E_l/E_{l-1}})$.
\end{enumerate}

\medskip
In the case where
$X$ is generically smooth and $(E, h)$ is
a $C^{\infty}$-hermitian locally free coherent sheaf on $X$,
for a non-zero $\OO_X$-subsheaf $G$ of $E$, we set
\[
\hat{\mu}_{\overline{H}}(G \hookrightarrow E, h) = \frac{\adeg(\acherncl_1(\overline{H}_1) 
\cdots \acherncl_1(\overline{H}_d) \cdot
\acherncl_1(G \hookrightarrow E, h))}{\rank G}.
\]

\bigskip
The purpose of this section is to prove the following unique existence
of an arithmetic Harder-Narasimham filtration:

\begin{Theorem}
\label{thm:existence:HN}
Let $X$ be a normal and projective arithmetic variety.
Let $(E, h)$ be a birationally $C^{\infty}$-hermitian torsion free
coherent sheaf on $X$. 
Let $\overline{H} = (\overline{H}_1, \ldots, \overline{H}_d)$ be a fine
sequence of nef $C^{\infty}$-hermitian invertible sheaves.
Then there exists uniquely
an arithmetic Harder-Narasimham filtration of $(E, h)$ with respect to $\overline{H}$.
Moreover, if $(E, h)$ is not
arithmetically $\mu$-semistable with respect to $\overline{H}$,
then a maximal destabilizing sheaf of $(E, h)$ is unique.
\end{Theorem}

We need several lemmas to prove the above theorem.

\begin{Lemma}
\label{lem:comp:birat:HN}
Let $(E, h)$ and $(E', h')$ be birationally 
$C^{\infty}$-hermitian torsion free
coherent sheaves on normal projective arithmetic varieties $X$ and
$X'$ respectively.
Let $\overline{H} = (\overline{H}_1, \ldots, \overline{H}_d)$ be a fine
sequence of nef $C^{\infty}$-hermitian invertible sheaves on $X$.
We assume that there is a birational morphism
$\pi : X' \to X$ and $(E, h)$ is dominated by
$(E', h')$ by means of $\pi : X' \to X$.
Then we have the followings:
\begin{enumerate}
\renewcommand{\labelenumi}{(\arabic{enumi})}
\item
$(E, h)$ is arithmetically $\mu$-semistable with respect to $\overline{H}$ if and only if
so is $(E', h')$ with respect to $\pi^*(\overline{H})$.

\item
Let $F$ be a saturated $\OO_X$-subsheaf of $E$ and $F'$
the corresponding saturated $\OO_{X'}$-subsheaf of $E'$.
Then $F$ is a maximal destabilizing sheaf of $(E, h)$ with respect
to $\overline{H}$ if and only if
so is $F'$ with respect to $\pi^*(\overline{H})$.

\item
Let $0 = E_0 \subsetneq E_1 \subsetneq \cdots \subsetneq E_l = E$
be a saturated filtration of $E$ and 
$0 = E'_0 \subsetneq E'_1 \subsetneq \cdots \subsetneq E'_l = E'$
the corresponding saturated filtration of $E'$.
Then $0 = E_0 \subsetneq E_1 \subsetneq \cdots \subsetneq E_l = E$
is a Harder-Narasimham filtration with respect to $\overline{H}$ if and only if
so is $0 = E'_0 \subsetneq E'_1 \subsetneq \cdots \subsetneq E'_l = E'$
with respect to $\pi^*(\overline{H})$.
\end{enumerate}
\end{Lemma}

\Proof
This is a consequence of Proposition~\ref{prop:birat:dominant}.
\QED

\begin{Lemma}
\label{lem:slope:comp}
Let $(E, h)$ be a birationally $C^{\infty}$-hermitian torsion free
coherent sheaf on a normal projective arithmetic variety $X$. 
If $(E, h)$ is not
arithmetically $\mu$-semistable with respect to $\overline{H}$ and
$F$ is a maximal slope sheaf of $(E, h)$, then
\[
\hat{\mu}_{\overline{H}}(F, h_{F\hookrightarrow E}) > 
\hat{\mu}_{\overline{H}}(E/F, h_{E \twoheadrightarrow E/F}).
\]
\end{Lemma}

\Proof
We set $a = \rank(F)$ and $b = \rank (E/F)$.
Then
\[
\hat{\mu}_{\overline{H}}(E, h)
= \frac{a}{a+b}\hat{\mu}_{\overline{H}}(F, h_{F\hookrightarrow E})
+ \frac{b}{a+ b}\hat{\mu}_{\overline{H}}(E/F, h_{E  \twoheadrightarrow E/F}).
\]
Thus, since
$\hat{\mu}_{\overline{H}}(F, h_{F\hookrightarrow E})
> \hat{\mu}_{\overline{H}}(E, h)$,  we get our lemma.
\QED

\begin{Lemma}
\label{lem:HN:birat}
Let $(E, h)$ be a birationally $C^{\infty}$-hermitian torsion free
coherent sheaf on a normal projective arithmetic variety $X$. 
Let $\overline{H} = (\overline{H}_1, \ldots, \overline{H}_d)$ be a fine
sequence of nef $C^{\infty}$-hermitian invertible sheaves.
Then there are a model $(E', h')$ of $(E, h)$ in terms of
a birational morphism $\mu : Y \to X$  of
normal projective arithmetic varieties and
a Harder-Narasimham filtration
\[
0 = E'_0 \subsetneq E'_1 \subsetneq \cdots \subsetneq E'_l = E'
\]
of $(E', h')$ with respect to $\mu^*(\overline{H})$ such that
$E'_i/E'_{i-1}$ is locally free for every $i = 1, \ldots, l$.
\end{Lemma}

\Proof
Let $(E', h')$ be a model of $(E, h)$ in terms of $\mu : Y \to X$.
By Proposition~\ref{prop:birat:dominant},
$(E, h)$ is arithmetically $\mu$-semistable with respect to $\overline{H}$
if and only if so is $(E', h')$
with respect to $\mu^*(\overline{H})$.
Thus we may assume that
$(E, h)$ is not arithmetically $\mu$-semistable with respect to $\overline{H}$.
Let $E'_1$ be a maximal destabilizing sheaf of $(E', h')$.
Considering Proposition~\ref{prop:birat:dominant} and
a suitable  birational morphism $\mu' : Y' \to Y$ of normal, projective and
generically smooth arithmetic varieties to remove the pinching points of $E'/E'_1$,
we may assume that $E'_1$ and $E'/E'_1$ are locally free.
If $(E'/E'_1, h'_{E' \twoheadrightarrow E'/E'_1})$
is arithmetically $\mu$-semistable, then we are done.
Otherwise, let $E'_2$ be a saturated $\OO_{Y}$-subsheaf of $E'$ such that
$E'_1 \subsetneq E'_2$ and $E'_2/E'_1$ is a maximal destabilizing sheaf of
$(E'/E'_1, h'_{E' \twoheadrightarrow E'/E'_1})$.
Changing $Y$ as before, we may assume that $E'_2$ and $E'/E'_2$ are locally free.
Moreover, by Lemma~\ref{lem:slope:comp},
\begin{multline*}
\hat{\mu}_{\mu^*(\overline{H})}(E'_1,  h_{E'_1 \hookrightarrow E'}) =
\hat{\mu}_{\mu^*(\overline{H})}(E'_1,  (h_{E'_2 \hookrightarrow E})_{E'_1 \hookrightarrow E'_2}) \\
> \hat{\mu}_{\mu^*(\overline{H})}(E'_2/E'_1, (h_{E'_2 \hookrightarrow E})_{E'_2 \twoheadrightarrow E'_2/E'_1}).
\end{multline*}
Thus, continuing this construction, we have our lemma.
\QED

\begin{Lemma}
\label{lem:max:slope:HN}
Let $(E, h)$ be a $C^{\infty}$-hermitian locally free
coherent sheaf on a normal projective and generically smooth arithmetic variety $X$. 
Let $\overline{H} = (\overline{H}_1, \ldots, \overline{H}_d)$ be a fine
sequence of nef $C^{\infty}$-hermitian invertible sheaves.
Let
$0 = E_0 \subsetneq E_1 \subsetneq \cdots \subsetneq E_l = E$ be an arithmetic
Harder-Narasimham filtration of $(E, h)$ such that
$E_i/E_{i-1}$ is locally free for every $i=1, \ldots, l$.
If $F$ is a maximal slope sheaf of $(E, h)$,
then $F \subseteq E_1$ and $\hat{\mu}_{\overline{H}}(F \hookrightarrow E, h) = 
\hat{\mu}_{\overline{H}}(E_1 \hookrightarrow E, h)$.
\end{Lemma}

\Proof
We choose $i$ such that $F \subseteq E_i$ and $F \not\subseteq E_{i-1}$.
We assume that $i \geq 2$.
Let $Q$ be the image of $F \to E_i/E_{i-1}$.
Let $h_Q$ be the quotient metric of $Q$ induced by $h_{F\hookrightarrow E}$ and $F \to Q$,
that is, $h_Q = (h_{F \hookrightarrow E})_{F  \twoheadrightarrow Q}$.
Then, by virtue of Lemma~\ref{lem:comp:inq:two:metrics},
\[
\hat{\mu}_{\overline{H}}(Q, h_Q) \leq
\hat{\mu}_{\overline{H}}(Q \hookrightarrow E_{i}/E_{i-1}, h_{E_{i}/E_{i-1}}).
\]
On the other hand, since $(F, h_{F\hookrightarrow E})$ 
and $(E_{i}/E_{i-1}, h_{E_{i}/E_{i-1}})$ are arithmetically
$\mu$-semistable,
\[
\hat{\mu}_{\overline{H}}(F, h_{F\hookrightarrow E}) \leq
\hat{\mu}_{\overline{H}}(Q, h_Q)
\]
and
\[
\hat{\mu}_{\overline{H}}(Q \hookrightarrow E_{i}/E_{i-1}, h_{E_{i}/E_{i-1}})
\leq \hat{\mu}_{\overline{H}}(E_{i}/E_{i-1}, h_{E_{i}/E_{i-1}}).
\]
Therefore,
\[
\hat{\mu}_{\overline{H}}(F, h_{F\hookrightarrow E}) \leq 
\hat{\mu}_{\overline{H}}(E_{i}/E_{i-1}, h_{E_{i}/E_{i-1}})
< \hat{\mu}_{\overline{H}}(E_1, h_{E_1 \hookrightarrow E}),
\]
which contradicts to the maximality of $\hat{\mu}_{\overline{H}}(F, h_{F\hookrightarrow E})$.
Thus $F \subseteq E_1$.
Moreover, since $(E_1, h_{E_1\hookrightarrow E})$ is arithmetically $\mu$-semistable,
$\hat{\mu}_{\overline{H}}(F, h_{F\hookrightarrow E}) \leq 
\hat{\mu}_{\overline{H}}(E_1, h_{E_1\hookrightarrow E})$.
Therefore $\hat{\mu}_{\overline{H}}(F, h_{F\hookrightarrow E}) = 
\hat{\mu}_{\overline{H}}(E_1, h_{E_1\hookrightarrow E})$
by  the maximality of $\hat{\mu}_{\overline{H}}(F, h_{F\hookrightarrow E})$.
\QED

\bigskip
Let us start the proof of Theorem~\ref{thm:existence:HN}.
The existence of a Harder-Narasimham filtration is
a consequence of Lemma~\ref{lem:HN:birat} and Proposition~\ref{prop:birat:dominant}.
Let us see the uniqueness of a Harder-Narasimham filtration.
Clearly we may assume that $(E, h)$ is not arithmetically $\mu$-semistable.
Let $0 = E_0 \subsetneq E_1 \subsetneq \cdots \subsetneq E_l = E$ and
$0 = G_0 \subsetneq G_1 \subsetneq \cdots \subsetneq G_{l'} = E$ be
Harder-Narasimham filtration of $(E, h)$.
Let $(E', h')$ be a model of $(E, h)$ in terms of $\mu : Y \to X$.
Let $0 = E'_0 \subsetneq E'_1 \subsetneq \cdots \subsetneq E'_l = E'$ and
$0 = G'_0 \subsetneq G'_1 \subsetneq \cdots \subsetneq G'_{l'} = E'$ be
corresponding Harder-Narasimham filtration of $(E', h')$ with
$0 = E_0 \subsetneq E_1 \subsetneq \cdots \subsetneq E_l = E$ and
$0 = G_0 \subsetneq G_1 \subsetneq \cdots \subsetneq G_{l'} = E$
respectively. By taking a birational morphism $\mu' : Y' \to Y$,
we may assume that $E'_i/E'_{i-1}$ and $G'_j/G'_{j-1}$ are
locally free for all $i = 1, \ldots, l$ and $j = 1, \ldots, l'$.
Let $F'$ be a maximal destabilizing sheaf of $(E', h')$.
Then, by Lemma~\ref{lem:max:slope:HN},
$F' \subseteq E'_1$ and $\hat{\mu}_{\mu^*(\overline{H})}(F', h_{F' \hookrightarrow E'}) = 
\hat{\mu}_{\mu^*(\overline{H})}(E'_1, h_{E'_1 \hookrightarrow E'})$.
Thus $F' = E'_1$. In  the same way, $F' = G'_1$.
Hence, by considering a Harder-Narasimham filtration of $(E'/F', h_{E' \twoheadrightarrow E'/F'})$
and induction on the rank, we have $l = l'$ and $E'_i = G'_i$ for all $i$.

The above observation also show the uniqueness of
a maximal destabilizing sheaf.
\QED

\end{document}